\magnification\magstephalf
\input amstex
\documentstyle{pub3}
\loadbold

\firstpage{369}
\lastpage{372}
\headtomus{50}
\headfasc{3-4}
\headyear{1997}

\headauthor {Victor Bovdi \es M. M. Parmenter}
\headtitle{Symmetric units in integral group rings}
\smalltitle{Symmetric units in integral group rings}

\topmatter
\title
Symmetric units in integral group rings
\endtitle
\author By {\smc VICTOR BOVDI} (Ny\'{\i}regyh\'aza) \es
{\smc M. M. PARMENTER} (St\. John's)
\endauthor
\received
Received October 14, 1996
\endreceived
\address
       Victor Bovdi\newline
       Department of Mathematics\newline
       Bessenyei Teachers College\newline
       4401 Ny\'{\i}regyh\'{a}za\newline
       Hungary
\endaddress
\email vbovdi\@math.klte.hu
\endemail
\address
       M. M. Parmenter\newline
       Department of Mathematics and Statistics\newline
       Memorial University of Newfoundland\newline
       St. John's, Newfoundland\newline
       Canada A1C 5S7
\endaddress
\subjclass
16534
\endsubjclass
\abstract
In this paper, we study the question of when the symmetric units in an integral
group ring $\Bbb Z\Bbb G$ form a multiplicative group. When $G$ is periodic,
necessary and sufficient conditions are given for this to occur.
\endabstract
\endtopmatter

\document

\footnotetext""{Research supported by the Hungarian National Foundation
for Scientific Research, Grant No. F015470, and by NSERC grant A8775, Canada.}

\heading {1. Introduction} \endheading

Let $U(KG)$ be the group of units of the group ring $KG$ of the group $G$ over
a commutative ring $K$.  The anti-automorphism $g \to g^{-1}$ of $G$ extends
linearly to an anti-automorphism $a \to a^*$ of $KG$. Let
$S_*(KG)=\{x\in U(KG)\mid x^*=x\}$ be the set of all symmetric units of
$U(KG)$.

\smallskip
The subgroup $U_*(KG) = \{x \in U(KG) \mid xx^* = 1\}$ is called the
{\it unitary\/} subgroup of $U(KG)$. It is easy to see ([4], Proposition 1.3)
that if $K=\Bbb{Z}$ then $U_* (\Bbb{Z} G)$ is trivial, i.e\.
$U_* (\Bbb{Z} G)=\pm G$. If $U(\Bbb{Z} G)\ne\pm G$, then in $U(\Bbb{Z} G)$
there always exist nontrivial symmetric units, for example $xx^*$ where $x$ is
a nontrivial unit in $U(\Bbb{Z} G)$.

\smallskip
In this paper we answer the question:  for which groups $G$ do the symmetric
units of the integral group ring $\Bbb{Z} G$ form a multiplicative group?
If $K$ is a commutative ring of characteristic $p$ and $G$ is a locally
finite $p$-group this question for $KG$ was described in [2].

\proclaim {Lemma {\rm(see [2])}}
Let $K$ be a commutative ring and $G$ be an arbitrary group. If $S_*(KG)$ is a
subgroup in $U(KG)$ then $S_*(KG)$ is abelian and normal in $U(KG)$.
\endproclaim

\proclaim {Theorem}
If $S_*(\Bbb{Z} G)$ is a subgroup in $U(\Bbb{Z} G)$, then the set $t(G)$ of
elements of $G$ of finite order is a subgroup in $G$, every subgroup of $t(G)$
is normal in $G$ and $t(G)$ is either abelian or a hamiltonian $2$-group.
Conversely, suppose that the group $G$ satisfies the above conditions and
$G/t(G)$ is a right ordered group. Then $S_*(\Bbb{Z} G)$ is a subgroup in
$U(\Bbb{Z} G)$.
\endproclaim

\heading {2. Proof of the theorem}\endheading

If the subgroup $t(G)$ of the group $G$ has the given properties and the
quotient group $G/t(G)$ is right ordered, then by Theorem 5.2 [1]
$$
V(\Bbb{Z} G) = G\cdot V (\Bbb{Z} t(G)). $$
Hence, every element $u\in S_*(\Bbb{Z}G)$ can be written as $bw$, where
$b$ is an element of $G$ and $w \in U(\Bbb{Z} t (G))$.  Suppose that $b$
is of infinite order and $w = \alpha_1 g_1 +\ldots+\alpha_s g_s$. Then
$bw =w^*b^{-1}$ and
$\text{Supp}(bwb)=\{bg_1 b,\dots,bg_s b\}=\{g^{-1}_1,\dots, g^{-1}_{s}\}$.
Thus $bg_1b = g^{-1}_{i}$ and $(bg_1)^2 = g^{-1}_{i} g_1$ is an element
of finite order, which is a contradiction.

\smallskip
We conclude that $S_*(\Bbb{Z} G) \subseteq U (\Bbb{Z} t(G))$.  If $t(G)$ is
abelian then $S_*(\Bbb{Z} G)$ is a subgroup.  On the other hand, if $t(G)$
is a hamiltonian 2-group then by Corollary 2.3 in [4], $V(\Bbb{Z} t(G)) = t(G)$
and so $S_*(\Bbb{Z} G)$ coincides with the centre of $t(G)$ and is again a
subgroup.

\smallskip
So now we assume that $S_*(\Bbb{Z} G)$ is a subgroup in $U(\Bbb{Z} G)$. We
first show that any subgroup of $t(G)$ is normal in $G$ (this also proves that
$t(G)$ is a subgroup of $G$). If not, then there exist $x\in t(G),y \in G$ with
$y^{-1} xy \notin\langle x\rangle$. But then $u=1+(1-x)y\hat{x}$ is a
nontrivial bicyclic unit in $\Bbb{Z}G$ (where $\hat{x}=1+x+\ldots+ x^{n-1}$,
$n=o(x))$, and {\smc Marciniak} and {\smc Sehgal} proved in [3] that
$\langle u,u^*\rangle$
is a nonabelian free subgroup of $U(\Bbb{Z} G)$. In particular, this means
that $uu^* \neq u^*u$ and that $uu^*,u^*u$ do not commute with each other.
Since $uu^*$ and $u^*u$ are in $S_*(\Bbb{Z} G)$, this contradicts the lemma.

\smallskip
We now have that $t(G)$ is either abelian or hamiltonian.  To finish the proof,
we need only to show that if $Q=\langle a,b\mid a^4=1,a^2=b^2,ba=a^3 b\rangle$
is the usual quaternion group and $g$ is of odd prime order $p$, then
$Q\times\langle g\rangle$ contains a pair of noncommuting symmetric
units.

\smallskip
Recall ([4], p\. 34) that if $x$ is of order $n$ in $G$ and $(i,n)=(j,n)= 1$,
and $ik \equiv 1\pmod n$, then
$$
u=(1+x^j+\ldots+x^{j(i-1)})(1+x^i+\ldots+x^{i(k-1)})+\frac{1-ik}n\hat{x} $$
is a (Hoechsmann) unit in $\Bbb{Z} G$.

\smallskip
First assume $p \neq 3$.  Then $ag$ and $bg$ are of order $4p$, and setting
$i = j = 3$ (and $3k \equiv 1 \pmod {4p}$) we obtain units
$$
\align
u & =(1+(ag)^3+(ag)^6)(1+(ag)^3+\ldots+(ag)^{3(k-1)})
 +\frac{1 - 3k}{4p} \widehat{ag}\\ \vspace {3pt}
v & =(1+(bg)^3+(bg)^6)(1+(bg)^3+\ldots+(bg)^{3(k-1)})
 +\frac{1 - 3k}{4p} \widehat{bg}\,.\endalign $$

Now $u_1 = (ag)^{-2} u$ and $v_1 = (bg)^{-2} v$ are symmetric units.  We
claim that $u_1$ and $v_1$ do not commute.  Since $(ag)^{-2}$ and $(bg)^{-2}$
are central, this is equivalent to showing that $u$ and $v$ do not commute.

\smallskip
Since $\frac{1-3k}{4p}\widehat{ag}$ and $\frac{1-3k}{4p}\widehat{bg}$
are central, this is equivalent to showing that $u_2$ and $v_2$ do not commute
where
$$
\align
u_2 & =(1+(ag)^3+(ag)^6)(1+(ag)^3+\dots+(ag)^{3(k-1)})\\ \vspace{3pt}
& =1+2(ag)^3+3(ag)^6+\ldots+3(ag)^{3(k-1)}+2(ag)^{3k}+(ag)^{3(k+1)}\\
 \vspace{5pt}
v_2 & =1+2(bg)^3+3(bg)^6+\ldots+3(bg)^{3(k-1)}+2(bg)^{3k}+(bg)^{3(k+1)}\,.
\endalign $$

Since all terms with even exponents are central, this is equivalent to showing
that $u_3$ and $v_3$ do not commute where
$$
\align
u_3 & =2 (ag)^3 + 3(ag)^9 +\ldots+ 3(ag)^{3(k-2)} + 2 (ag)^{3k}\\ \vspace {3pt}
v_3 & =2 (bg)^3 + 3 (bg)^9 +\ldots+3 (bg)^{3(k-2)} + 2 (bg)^{3k}.\endalign $$

But in $u_3v_3$ only 4 products are not divisible by 3. Since
$3k \equiv 1\pmod{4p}$, these reduce to $4abg^6+8a^3bg^4+4abg^2$. In $v_3u_3$,
the same products reduce to $4a^3bg^6+8abg^4+4a^3bg^2$. Because all other
products are divisible by 3, we see $u_3v_3 \neq v_3u_3$.

\smallskip
If $p = 3$, the same argument works with $i = j = k = 5$.  In this case, direct
calculation shows that if $u$ and $v$ are defined as before, the symmetric
units $(ag)^4 u$ and $(bg)^4 v$ do not commute.
\hfill$\qed$

\medskip
Note that when $G$ is periodic, the theorem shows that $S_*(\Bbb{Z} G)$ is a
subgroup only in the obvious cases -- namely when $G$ is either abelian or
a hamiltonian 2-group.

\smallskip
We remark that it is possible to avoid using the result from [3] and to prove
that every subgroup of $t(G)$ is normal in $G$ by a direct argument instead. We
have decided to use [3] in order to indicate how useful the Marciniak--Sehgal
result can be.

\enddocument

\enddocument